
\input amstex.tex
\documentstyle{amsppt}
\magnification=\magstep1
\hsize=12.5cm
\vsize=18cm
\hoffset=1cm
\voffset=2cm
\def\DJ{\leavevmode\setbox0=\hbox{D}\kern0pt\rlap
 {\kern.04em\raise.188\ht0\hbox{-}}D}
\footline={\hss{\vbox to 2cm{\vfil\hbox{\rm\folio}}}\hss}
\nopagenumbers
\font\ff=cmr8
\def\txt#1{{\textstyle{#1}}}
\baselineskip=13pt
\def\hf{{\textstyle{1\over2}}}
\def\a{\alpha}
\def\d{{\,\roman d}}
\def\e{\varepsilon}
\def\f{\varphi}
\def\G{\Gamma}
\def\k{\kappa}
\def\s{\sigma}

\def\z{\zeta}
\def\={\;=\;}
\def\zx{\zeta(\hf+ix)}
\def\zt{\zeta(\hf+it)}

\def\no{\noindent}  
\def\R{\Re{\roman e}\,}  \def\s{\sigma}
\def\Z{{\Cal Z}}

\def\no{\noindent}

\font\teneufm=eufm10
\font\seveneufm=eufm7
\font\fiveeufm=eufm5
\newfam\eufmfam
\textfont\eufmfam=\teneufm
\scriptfont\eufmfam=\seveneufm
\scriptscriptfont\eufmfam=\fiveeufm
\def\mathfrak#1{{\fam\eufmfam\relax#1}}

\font\tenmsb=msbm10
\font\sevenmsb=msbm7
\font\fivemsb=msbm5
\newfam\msbfam
\textfont\msbfam=\tenmsb
\scriptfont\msbfam=\sevenmsb
\scriptscriptfont\msbfam=\fivemsb
\def\Bbb#1{{\fam\msbfam #1}}

\def \NN {\Bbb N}
\def \CC {\Bbb C}

\def \ZZ {\Bbb Z}

\def\rightheadline{{\hfil{\ff
The modified Mellin transform of powers of the zeta-function}\hfil\tenrm\folio}}

\def\leftheadline{{\tenrm\folio\hfil{\ff
Aleksandar Ivi\'c }\hfil}}
\def\emptyheadline{\hfil}
\headline{\ifnum\pageno=1 \emptyheadline\else
\ifodd\pageno \rightheadline \else \leftheadline\fi\fi}

\topmatter
\title THE MODIFIED MELLIN TRANSFORM OF POWERS
OF THE ZETA-FUNCTION \endtitle
\author   Aleksandar Ivi\'c
\bigskip
\endauthor
\dedicatory
Dedicated to the memory of Yu. V. Linnik
\enddedicatory
\address{
Aleksandar Ivi\'c, Katedra Matematike RGF-a
Universiteta u Beogradu, \DJ u\v sina 7, 11000 Beograd,
Serbia and Montenegro.
}
\endaddress
\email {\tt aivic\@ivic.matf.bg.ac.yu, ivic\@rgf.bg.ac.yu} \endemail
\abstract
{The modified Mellin transform  $\Z_k(s) = \int_1^\infty |\zx|^{2k}x^{-s}\d x
\;(k\in\NN)$ is investigated. Analytic continuation and mean square estimates
of $\Z_k(s)$ are discussed, as well as
connections with power moments of $|\zx|$, with the special
emphasis on the cases $k = 1,2$.}
\endabstract
\endtopmatter

\heading
{\bf 1. Introduction}
\endheading

Integral transforms (Laplace, Fourier, Mellin among others)
play an important r\^ole in analytic number
theory. Of special interest in the theory of the Riemann
zeta-function $\zeta(s)$ are the Laplace transforms
$$
L_k(s) \;:=\; \int_0^\infty |\zx|^{2k}{\roman e}^{-sx}\d x
\qquad(k \in \NN,\,\s = \R s > 0)
\leqno(1.1)
$$
and the (modified) Mellin transforms
$$
\Z_k(s) \;:=\; \int_1^\infty |\zx|^{2k}x^{-s}\d x
\qquad(k \in \NN,\,\s = \R s \ge c(k) > 1),
\leqno(1.2)
$$
where $c(k)$ is such a constant for which the integral in (1.2) converges
absolutely. The term ``modified" Mellin transform seems appropriate, since
customarily the Mellin transform of $f(x)$ is defined as
$$
 F(s) := \int_0^\infty f(x)x^{s-1}\d x\qquad(s = \s + it\in\CC).\leqno(1.3)
$$
Note that the lower bound of integration in (1.2) is not zero,
as it is in (1.3). The choice
of unity as the lower bound of integration dispenses with convergence
problems at that point, while the appearance of the factor $x^{-s}$
instead of the customary $x^{s-1}$ is technically more convenient. Also
it may be compared with the discrete representation
$$
\zeta^{2k}(s) \= \sum_{n=1}^\infty d_{2k}(n)n^{-s}\qquad(\s >
1,\,k\in\NN),
$$
where $d_m(n)$ is the number of ways $n$ may be written as a product of
$m$  factors; $d(n) \equiv d_2(n)$ is the number of divisors of $n$.
Since we have (see   [6, Chapter 8])
$$
\int_0^T|\zt|^{2k}\d t \ll T^{(k+2)/4}\log^{C(k)}T\qquad(2 \le k \le 6),
\leqno(1.4)
$$
it follows that the integral defining ${\Cal Z}_k(s)$ is absolutely
convergent
for $\s > 1$ if $0 \le k \le 2$ and for $\s > (k+2)/4$ if $2 \le k \le 6$.

\medskip

E.C. Titchmarsh's well-known monograph [26, Chapter 7] gives a discussion
of $L_k(s)$ when $s = \s$ is real and $\s \to 0+$, especially detailed
in the cases $k=1$ and $k=2$. Indeed, a classical result of H. Kober
[19] says that, as $\s \to 0+$,
$$
L_1(2\s) = {\gamma-\log(4\pi\s)\over2\sin\s} +
 \sum_{n=0}^Nc_n\s^n + O(\s^{N+1})
\leqno(1.5)
$$
for any given integer $N \ge 1$, where the $c_n$'s are effectively
computable constants and $\gamma = 0.577\ldots\,$ is Euler's constant.
For complex values of $s$ the function $L_1(s)$ was studied by
F.V. Atkinson [1], and more recently by M. Jutila [17].
Atkinson [2] obtained the asymptotic formula, as $\s \to 0+$,
$$
L_2(\s) = {1\over\s}\left(A\log^4{1\over\s} + B\log^3{1\over\s}
+ C\log^2{1\over\s} + D\log {1\over\s} + E\right) + \lambda_2(\s),
\leqno(1.6)
$$
where
$
A = {1\over2\pi^2},\,B ={1\over\pi^{2}}
\Bigl(2\log(2\pi) - 6\gamma + {24\zeta'(2)\over\pi^{2}}\Bigr),
$
and
$
 \lambda_2(\s) \;\ll_\e\;\left({1\over\s}\right)^{{13\over14}+\e}.
$
We note that here and later $\e$ denotes arbitrarily small constants,
not necessarily the same ones at each occurrence.

The author [8] gave explicit, albeit complicated expressions for
the remaining coefficients $C,D$ and $E$ in (1.6). More importantly,
he applied a result on the fourth moment of $|\zt|$, obtained jointly
with Y. Motohashi (see the second bound in (6.2)), to establish that
$$
\lambda_2(\s) \;\ll\; \s^{-1/2}\qquad(\s\to 0+).
$$

\medskip
For $k \ge 3$ not much is known about $L_k(s)$, even when $s = \s \to 0+$.
This is not surprising, since not much is known about upper bounds for
$$
I_k(T) \;:=\; \int_0^T|\zt|^{2k}\d t\qquad(k \in \NN), \leqno(1.7)
$$
where $k \ge 3\,$. For a discussion on $I_k(T)$ the reader is referred
to the author's monographs [6] and [7]. One trivially has
$$
I_k(T) \le {\roman e}\int_0^\infty|\zt|^{2k}
{\roman e}^{-t/T}\d t = {\roman e}L_k\left({1\over T}\right).
\leqno(1.8)
$$
Thus any nontrivial bound of the form
$$
L_k(\s) \;\ll_\e\; \left({1\over\s}\right)^{c_k+\e}\qquad(\s\to 0+,\,
c_k \ge 1)\leqno(1.9)
$$
gives, in view of (1.8) ($\s = 1/T$), the bound
$$
I_k(T) \;\ll_\e\; T^{c_k+\e}.\leqno(1.10)
$$
Conversely, if (1.10) holds, then we obtain (1.9) from the identity
$$
L_k\left({1\over T}\right) \= {1\over T}\,\int_0^\infty I_k(t)
{\roman e}^{-t/T}\d t,
$$
which is easily established by integration by parts.

Note that, by the change of variable $x = {\roman e}^t,\,z = s-1$, (1.2) becomes
$$
\int_0^\infty |\z(\hf + i{\roman e}^t)|^{2k}{\roman e}^{-zt}\d t\qquad
(\R z > 0),
$$
which is the Laplace transform of $|\z(\hf + i{\roman e}^t)|^{2k}$. Indeed, it
is well-known that the Laplace and Mellin transforms are closely connected, as
(by a change of variable) both of them can be regarded as special cases of
Fourier transforms, and their theory built from the theory of Fourier transforms.

\medskip
The aim of this paper is to give an account on the known results for $\Z_k(s)$
and to prove some new results. The function $\Z_1(s)$ was investigated by
Ivi\'c, Jutila and Motohashi [16], and later by M. Jutila [17] and the author [13].
The function $\Z_2(s)$ was introduced by Y. Motohashi [23] (see also his monograph
[24]). It was later investigated in [16], as well as in the author's works [11]
and [12]. The papers [11] and [16] also contain material
on the general function $\Z_k(s)$.

\bigskip\medskip
\heading {\bf  2.
The analytic continuation of $\Z_k(s)$ }
  \endheading
Remarks on the general problem of analytic continuation of $\Z_k(s)$
were given in [11] and [16]. We start here by proving a general result, which
links the problem to the moments of $|\zt|$. This is

\bigskip
THEOREM 1. {\it Let $k\in\NN$ be fixed. The bound
$$
\int_0^T|\zt|^{2k}\d t \;\ll_\e\; T^{c+\e}\leqno(2.1)
$$
holds for some constant $c$, if and only if $\Z_k(s)$ is regular for
$\R s > c$, and for any given $\e>0$}
$$
\Z_k(c+\e+it) \;\ll_\e\; 1.\leqno(2.2)
$$

\bigskip
{\bf Proof of Theorem 1}. The constant $c$ must satisfy $c\ge1$ in view of
the known lower bounds for moments of $|\zt|$ (see e.g., [6, Chapter 9]).
Suppose that (2.1) holds. Then we have
$$
\int_X^{2X}|\zx|^{2k}x^{-s}\d x \;\ll
\; X^{-\s}\int_0^{2X}|\zx|^{2k}\d x
\ll_\e X^{c-\s+\e/2},
$$
where $\s = \R s$. Therefore
$$
\eqalign{\int_1^\infty |\zx|^{2k}x^{-s}\d x& \;=\;
\sum_{j=0}^\infty \,\int_{2^j}^{2^{j+1}}|\zx|^{2k}x^{-s}\d x\cr&
\;\ll\; \sum_{j=0}^\infty 2^{j(c-\s+\e/2)} \;\ll_\e\; 1
\cr}
$$
if $\s = c+\e$, since $\sum_j2^{-\e j/2}$ converges.
This shows that $\Z_k(s)$ is regular for $\s > c$ and that
(2.2) holds.

Conversely, suppose that $\Z_k(s)$ is regular for $\s > c$ and that
(2.2) holds. Using the classical integral ($\int_{(d)}$ denotes
integration over the line $\R s = d$)
$$
{\roman e}^{-x} = {1\over2\pi i}\int_{(d)}x^{-s}\G(s)\d s\qquad
(\R x > 0,\,d>0),
$$
we have
$$
\eqalign{
\int_1^\infty {\roman e}^{-x/T}|\zx|^{2k}\d x &=
\int_1^\infty {1\over2\pi i}\int_{(c+\e)}\G(s)\bigl({x\over T}\bigr)^{-s}\d s
|\zx|^{2k}\d x\cr&
= {1\over2\pi i}\int_{(c+\e)}\G(s) T^s\Z_k(s)\d s \ll_\e T^{c+\e}\cr},
$$
by absolute convergence and the fast decay of the gamma-function. This
yields
$$
\eqalign{&
\int_0^T|\zx|^{2k}\d x \le O(1) + {\roman e}\int_1^T
{\roman e}^{-x/T}|\zx|^{2k}\d x\cr&
\ll 1 + \int_1^\infty {\roman e}^{-x/T}|\zx|^{2k}\d x
\ll_\e T^{c+\e},
\cr}
$$
which proves (2.1).

\bigskip
{\bf Corollary 1}. The Lindel\"of hypothesis ($|\zt| \ll_\e |t|^\e$)
is equivalent to the statement that, for every $k\in\NN$, $\Z_k(s)$ is
regular for $\s>1$ and satisfies $\Z_k(1+\e+it) \ll_{k,\e} 1$.

\medskip
Indeed, the Lindel\"of hypothesis is equivalent (see e.g., [6, Section 1.9])
to (2.1) with $c=1$ for every $k\in\NN$. Therefore the assertion follows
from Theorem 1.

\bigskip
{\bf Corollary 2}. If we define
$$
\s_k := \inf\left\{\,d_k \;:\; \int_0^T|\zt|^{2k}\d t \ll_{k}
T^{d_k}\,\right\},\leqno(2.3)
$$
$$
\rho_k := \inf\left\{\,r_k \;:\; \Z_k(s) \;\text {is regular for}\;
\R s > r_k\,\right\},\leqno(2.4)
$$
then
$$
\rho_k \;=\;\s_k,\quad \s_k \ge\;1.\leqno(2.5)
$$

\medskip
Note that from the classical bound
$\int_0^T|\zt|^2\d t \ll T\log T$ and (1.4) we obtain
$$
\s_1 = \s_2 = 1,\quad \s_k \le {k+2\over4}\quad(3\le k \le 6),\leqno(2.6)
$$
and upper bounds for $\s_k$ when $k>7$ may be obtained by using results
on the corresponding power moments of $|\zt|$ (see [6, Chapter 8]).
The Lindel\"of hypothesis may be reformulated as $\s_k = 1\;(\forall k \ge 1)$.

\smallskip
Thus at present we have two situations regarding analytic continuation
of $\Z_k(s)$:

\smallskip
a) For $k = 1,2,$ one can obtain analytic continuation of $\Z_k(s)$ to
the left of $\R s = 1$ (in fact to $\CC$). This will be discussed in
Section 5 and Section 6, respectively.

\smallskip

b) For $k>2$ only upper bounds for $\s_k$ (cf. (2.6)) are known. A
challenging problem is to improve these bounds, which would entail
progress on bounds of power moments of $|\zt|$, one of the central
topics in the theory of $\z(s)$.

\medskip
In what concerns power moments of $|\zt|$ one expects,
for any fixed $k \in\NN$,
$$
\int_0^T|\zt|^{2k}\d t = TP_{k^2}(\log T) + E_k(T)\leqno(2.7)
$$
to hold, where it is generally assumed that
$$
P_{k^2}(y) = \sum_{j=0}^{k^2}a_{j,k}y^j\leqno(2.8)
$$
is a polynomial in $y$ of degree $k^2$ (the integral in (2.7) is
$\gg_k T\log^{k^2}T$; see e.g., [6, Chapter 9]). The function
$E_k(T)$ is to be considered as the error term in (2.7), namely
one supposes that
$$
E_k(T) \= o(T)\qquad(T \to \infty).\leqno(2.9)
$$
So far (2.7)--(2.9) are known to hold only for $k = 1$ and $k = 2$
(see [7] and [24] for a comprehensive account). Therefore in view of the
existing knowledge on the higher moments of $|\zt|$,
embodied in (1.4), at present
the really important cases of (2.7) are $k = 1$ and $k = 2$.

\medskip
In case (2.7)--(2.9) hold, this may be used to obtain the
analytic continuation of $\Z_k(s)$ to the region $\s\ge 1$ (at least).
Indeed, by using (2.7)-(2.9) we have
$$
\eqalign{&
\Z_k(s) = \int_1^\infty |\zx|^{2k}x^{-s}\d x = \int_1^\infty
x^{-s}\d\left(xP_{k^2}(\log x) + E_k(x)\right)\cr&
= \int_1^\infty (P_{k^2}(\log x) + P'_{k^2}(\log x))x^{-s}\d x
- E_k(1) + s\int_1^\infty E_k(x)x^{-s-1}\d x.\cr}\leqno(2.10)
$$
But for $\R s > 1$ change of variable $\log x = t$ gives
$$
\eqalign{&
\int_1^\infty (P_{k^2}(\log x) + P'_{k^2}(\log x))x^{-s}\d x\cr&
 = \int_1^\infty \left\{\sum_{j=0}^{k^2}a_{j,k}\log^jx +
\sum_{j=0}^{k^2-1}(j+1)a_{j+1,k}\log^jx\right\}x^{-s}\d x\cr&
= \int_0^\infty \left\{\sum_{j=0}^{k^2}a_{j,k}t^j +
\sum_{j=0}^{k^2-1}(j+1)a_{j+1,k}t^j\right\}{\roman e}^{-(s-1)t}\d t\cr&
= {a_{k^2,k}(k^2)!\over(s-1)^{k^2+1}} +
\sum_{j=0}^{k^2-1}(a_{j,k}j! + a_{j+1,k}(j+1)!)(s-1)^{-j-1}.
\cr}\leqno(2.11)
$$
Hence inserting (2.11) in (2.10) and using (2.9)
we obtain  the analytic continuation of $\Z_k(s)$
to the region $\s\ge1$. As we know (see [7] and [24]) that
$$
\int_1^T E_1^2(t)\d t \ll T^{3/2},\qquad
\int_1^T E_2^2(t)\d t \ll T^{2}\log^{22}T,\leqno(2.12)
$$
it follows on applying the Cauchy--Schwarz inequality to the last
integral in (2.10) that (2.9)-(2.11) actually provides
the analytic continuation of $\Z_1(s)$ to the region $\R s > 1/4$, and of
$\Z_2(s)$ to $\R s > 1/2$.
\heading
{\bf 3. Recurrence relations and identities}
\endheading

There is a possibility to obtain analytic continuation
of $\Z_k(s)$ by using a recurrent relation involving
$\Z_r(s)$ with $r < k$, which was mentioned in [11] and [16].
This is

\medskip
THEOREM 2. {\it For $k\ge 2,\,r = 1,\ldots, k-1$, $\R s$ and $c = c(k,r)$ sufficiently
large, we have}
$$
{\Cal Z}_k(s) = {1\over2\pi i}\int_{(c)}{\Cal Z}_{k-r}(w){\Cal Z}_r(1+s-w)
\d w.\leqno(3.1)
$$

\medskip
{\bf Proof of Theorem 2}. For $\R(1-s)$ sufficiently large we have
$$
{\Cal Z}_k(1-s) = \int_0^\infty \zeta^*(x)x^{s-1}\d x,
$$
where $\zeta^*(x) = |\zx|^{2k}$ if $x \ge 1$ and zero otherwise.
Consequently Mellin inversion (see e.g., the Appendix of [6]) will give
$$
|\zx|^{2k} = {1\over2\pi i}\int_{(c)}{\Cal Z}_k(1-s)x^{-s}\d s,
\qquad(c \le c_0(k) < 0,\;x \ge 1).\leqno(3.2)
$$
Therefore, for $\,k,r \in \NN,\,k\ge2,\,1 \le r < k$, we obtain
$$\eqalign{&
\int_1^\infty |\zx|^{2k}x^{-s}\d x = \int_1^\infty |\zx|^{2r}
|\zx|^{2(k-r)}x^{-s}\d x\cr&
= \int_1^\infty |\zx|^{2r}\left({1\over2\pi i}\int_{(c)}{\Cal Z}_{k-r}
(1-w)x^{-w}\d w\right)x^{-s}\d x\cr&
= {1\over2\pi i}\int_{(c)}{\Cal Z}_{r}(w+s){\Cal Z}_{k-r}(1-w)\d w
\quad(\s \ge \s_0(k) > 1-c).\cr}
$$
Changing $1-w$ to $w$ we obtain (3.1).

\medskip
In particular, by using (1.4), we obtain the identities
$$\eqalign{
\Z_3(s) &\= {1\over2\pi i}\int_{(1+\e)}{\Z}_1(w)
\Z_2(1 - w + s)\d w\qquad(\s > {\txt{5\over4}}),\cr
{\Z}_4(s) &\= {1\over2\pi i}\int_{({5\over4}+\e)}{\Z}_2(w)
{\Z}_2(1 - w + s)\d w\qquad(\s > {\txt{3\over2}}).\cr}
$$

\medskip

The following result provides an integral representation for $\Z_k^2(s)$.
This is

\bigskip
THEOREM 3. {\it In the region of absolute convergence we have}
$$
\Z_k^2(s) = 2\int_1^\infty x^{-s}\left(\int_{\sqrt{x}}^x|\zx|^{2k}
\Bigl|\z\Bigl(\hf +i\,{x\over u}\Bigr)\Bigr|^{2k}{\d u\over u}\right)\d x.
\leqno(3.3)
$$

\bigskip
{\bf Proof of Theorem 3}. Set $f(x) = |\zx|^{2k}$ and make the change
of variables $xy = X,\,x/y = Y$, so that the absolute value of the
Jacobian of the transformation is equal to $1/(2Y)$. Therefore
$$\eqalign{
\Z_k^2(s) &= \int_1^\infty\int_1^\infty (xy)^{-s}f(x)f(y)\d x\d y\cr&
= {1\over2}\int_1^\infty X^{-s}\int_{1/X}^X{1\over Y}
f(\sqrt{XY}\,)f(\sqrt{X/Y}\,)\d Y\d X.\cr}
$$
But as we have ($y = 1/u$)
$$
\int_{1/x}^x
f(\sqrt{xy}\,)f(\sqrt{x/y}\,){\d y\over y} = \int_1^xf(\sqrt{x/u})
f(\sqrt{xu}\,){\d u\over u},
$$
we obtain that, in the region of absolute convergence, the identity
$$
\Z_k^2(s) = \int_1^\infty x^{-s}\left(\int_1^x f(\sqrt{xy}\,)
f(\sqrt{x/y}\,){\d y\over y}\right)\d x
$$
is valid. The inner integral here becomes, after the change of
variable $\sqrt{xy} = u$,
$$
2\int_{\sqrt{x}}^x f(u)f\bigl({x\over u}\bigr){\d u\over u},
$$
and (3.3) follows. The argument also shows that, for  $0 < a < b$
and any integrable function $f$ on $[a,\,b]$,
$$\eqalign{&
\left(\int_a^b f(x)x^{-s}\d x\right)^2\cr&
= 2\int_{a^2}^{b^2}x^{-s}\left\{\int_{\sqrt{x}}^{\min(x/a,b)}f(u)
f\bigl({x\over u}\bigr){\d u\over u}\right\}\d x.
\cr}
\leqno(3.4)
$$

\heading
{\bf 4. Multiple Dirichlet series}
\endheading

In the recent work [5] of Diaconu, Goldfeld and Hoffstein
the theory of multiple Dirichlet series is developed.
In particular, they consider the series
$$
Z(s_1,\cdots,s_{2m},w)
= \int_1^\infty \zeta(s_1+it)\cdots \zeta(s_m+it)\zeta(s_{m+1}-it)
\cdots\zeta(s_{2m}-it)t^{-w}\,{\roman d}t\leqno(4.1)
$$
connected with the Riemann zeta-function. Analytic properties of
this function, closely  connected to our function $\Z_k(s)$,
may be put to advantage to deal with the important problem of the
analytic continuation of the function $\Z_k(s)$ itself.
It is shown in [5] that (4.1) has meromorphic
continuation (as a function of $2m+1$ complex variables) slightly
beyond the region of absolute convergence, with a polar divisor at
$w=1$. It is also shown that (4.1) satisfies certain quasi-functional
equations, which are used to obtain meromorphic continuation to an
even larger region. Under the assumption that
$$
Z(\hf,\cdots,\hf,w) \;\equiv\; \Z_m(w)
$$
has holomorphic continuation to
the region $\Re{\roman e}\, w \ge 1$ (except for the pole at $w=1$
of order $m^2+1$),
the authors derive the conjecture
on the moments of the zeta-function on the critical line in the form
$$
\int_0^T|\zt|^{2k}\d t = (c_k + o(1))T\log^{k^2}T\qquad(T\to\infty),\leqno(4.2)
$$
where $k\ge2$ is a fixed integer and
$$
c_k = {a_kg_k\over\G(1+k^2)},\quad a_k = \prod_{p}(1-1/p)^{k^2}
\sum_{j=0}^\infty d^2_k(p^j),\quad g_k = (k^2)!\prod_{j=0}^{k-1}{j!\over(j+k)!}.
\leqno(4.3)
$$
The formulas (4.2)-(4.3) coincide with the conjecture from Random
Matrix Theory (see e.g., Keating--Snaith [18]).  This is a weak form
of (2.7)--(2.9), but the point here is in the explicit form of the
constant $c_k$ (it equals $a_{k^2,k}$ in (2.8)).

 In [5] the conjectural formula (4.3), in the final step,
 is derived from a Tauberian theorem. The Tauberian theorem is the following

\bigskip
PROPOSITION 1. {\it Let $F(x) \;(x\ge 1)$ be a non-decreasing continuous
function, and
$$
f(s) := \int_1^\infty F(u)u^{-s-1}\d u.
$$
Let
$$
P(s) = \gamma_M + \gamma_{M-1}(s-1) + \cdots + \gamma_0(s-1)^M\qquad
(\gamma_M \ne 0),
$$
and suppose that $f(s) - P(s)(s-1)^{-M-1}$ is holomorphic for
$\R s > 1$ and continuous for $\R s = 1$. Then}
$$
F(x) \;\sim\;{\gamma_M\over M!}x(\log x)^M\qquad(x\to\infty).\leqno(4.4)
$$

\bigskip
When $M=0$ this is the classical Wiener--Ikehara Tauberian theorem
(see e.g., J. Korevaar [20, Theorem III.4.1]). Note that (op. cit., eq. (4.1))
$\d s(v) = \d(s(v)-B)$ for  any constant $B$, hence we may in fact assume
that $F(x)$ above is actually non-negative.
In [5] there is no proof of Proposition 1 (it is attributed to
H. Stark's unpublished notes), so a discussion here seems in place.
J. Korevaar kindly pointed out to
me that M.A. Subhankulov [25] derives the general case under some
more stringent conditions than the ones given above (but quite
sufficient for applications involving the moments of $|\zt|$).
A generalized version of the Wiener--Ikehara Tauberian theorem
was obtained long ago by H. Delange [4], but it is not obvious
whether his arguments can be modified to yield the above result.
Prof. Korevaar also pointed out to me how the
general case, stated above, can be reduced to follow from his
proof of the Wiener--Ikehara Tauberian theorem, given in [20].
First we set $S(x) = F({\roman e}^x)$, and then we may assume that
$S(x)$ is non-decreasing, non-negative and $S(x) = 0$ for $x\le 1$. Instead
of (4.4) it is sufficient to prove that
$$
S(x) \;\sim\; {\gamma_M\over M!}{\roman e}^xx^M\qquad(x\to\infty).\leqno(4.5)
$$
Change of variable $u = {\roman e}^x,\, s = z+1$ gives
$$
f(s) = \int_1^\infty F(u)u^{-s-1}\d u
= \int_1^\infty S(x){\roman e}^{-x}\cdot {\roman e}^{-zx}\d x = \f(z) ,
$$
say, so that $\f(z)$ is the Laplace transform of $S(x){\roman e}^{-x}$. We have
that
$$
g(z) := \f(z) - {\gamma_M\over z^{M+1}} - \cdots - {\gamma_0\over z}\leqno(4.6)
$$
is regular for $x = \R z > 0$ and continuous for $x\to0+$. Integration of (4.6)
shows that the same holds for
$$\eqalign{
g_1(z) &:= \int_1^z g(w)\d w = -\int_1^\infty S(x)
{\roman e}^{-x}x^{-1}({\roman e}^{-zx}-{\roman e}^{-x})\d x\cr&
+ {\gamma_M\over Mz^M} + \cdots - \gamma_0\log z + C,\cr}
$$
where $C$ is a constant.
Continuing the process with $g_n(z) := \int_1^z g_{n-1}(w)\d w$, it follows that
the same holds for
$$
g_M(z) = (-1)^M\int_1^\infty S(x){\roman e}^{-x}x^{-M}{\roman e}^{-zx}\d x
- (-1)^M{\gamma_M\over M!z} + \G_M(z),
$$
say, where $\G_M(z)$ is also regular for $x>0$ and continuous for $x\to0+$.
Now it remains to follow the proof of the Wiener--Ikehara Tauberian theorem
[20], which follows from [17, Ch. 3, Proposition 4.3].
The change is that, instead of
$\s(x) = S(x){\roman e}^{-x}$, now we shall work with
$$
\s(x) = S(x){\roman e}^{-x}x^{-M}.
$$
By the non-decreasing property of $S(x)$ it follows that
$$
\s(u - v/\lambda) \ge \s(u - a/\lambda){\roman e}^{-2a/\lambda}
\left(u - a/\lambda\over u + a/\lambda\right)^M\leqno(4.7)
$$
for $u \ge u_0 \,(>0)$, $|v| \le a$ and $a,\lambda > 0$. Therefore by using (4.7)
we obtain
$$
\eqalign{
A &= \lim_{u\to\infty}\int_{-\infty}^{\lambda u}\s(u - v/\lambda)K(v)\d v
\ge \limsup_{u\to\infty}\int_{-a}^{a}\s(u - v/\lambda)K(v)\d v\cr&
\ge \limsup_{u\to\infty}\,\s(u - a/\lambda){\roman e}^{-2a/\lambda}
\left(u - a/\lambda\over u + a/\lambda\right)^M\int_{-a}^a K(v)\d v,\cr}
$$
where $K$ is the Fej\'er kernel
$K_\lambda(t) = \lambda K(\lambda t) = {\lambda\over2\pi}
\left({\sin \lambda t/2\over\lambda t/2}\right)^2$. This yields
$$
\limsup_{u\to\infty}\,\s(u) = \limsup_{u\to\infty}\,\s(u - a/\lambda)
\left(u - a/\lambda\over u + a/\lambda\right)^M \le {{\roman e}^{2a/\lambda}A
\over \int_{-a}^a K(v)\d v}.\leqno(4.8)
$$
From (4.8) it follows, taking e.g.,
$a = \sqrt{\lambda},\,\lambda\to\infty$, that
$$
\limsup_{u\to\infty}\,\s(u)  \;\le \;A,
$$
and $\liminf_{u\to\infty}\s(u)\ge A$ follows analogously as in Korevaar's book.
This proves that
$$
\lim_{u\to\infty}\,\s(u) = \lim_{u\to\infty} S(u){\roman e}^{-u}u^{-M}
 = A,
$$
hence (4.5) follows (since $A = \gamma_M/ M!$ in our case)
and Proposition 1 is proved. Another variant of proof,
involving also repeated integration, was kindly pointed out to me
by Prof. Y. Motohashi.

\medskip
Note that (4.4) of Proposition 1 gives the value of $c_k = a_{k^2,k}$,
the (conjectural) leading coefficient of the polynomial $P_{k^2}(y)$
in (2.7), but not the remaining coefficients.  Although a plausible
conjecture for the remaining coefficients has been recently given by
Random Matrix Theory (see J.B. Conrey et al. [3]), there seems to exist
no Tauberian theorem strong enough which would furnish these coefficients
from the knowledge of the analytic behaviour of $\Z_k(s)$ near $s=1$.

\heading
{\bf 5. The  function $\Z_1(s)$}
\endheading
We begin with the analytic continuation of $\Z_1(s)$. This is contained in

\bigskip
{THEOREM 4}. {\it The function $\Z_1(s)$
continues meromorphically to $\CC$, having only a double pole
at $s=1$, and at most simple poles at $s = -1,-3,\ldots\,$. The principal
part of its Laurent expansion at $s=1$ is given by
$$
{1\over(s-1)^2} + {2\gamma - \log(2\pi)\over s - 1},\leqno(5.1)
$$
where $\gamma = -\G'(1) = 0.577215\ldots\,$ is Euler's constant.}
\bigskip
{\bf Proof of Theorem 4}. It was shown in [16] that $\Z_1(s)$
continues analytically to a function that is regular for $\s > -3/4$. In [17]
M. Jutila proved
that $\Z_1(s)$ continues meromorphically to $\CC$, having only a double pole
at $s=1$ and at most double poles for $s = -1,-2,\ldots\,$. The present form
of Theorem 4 was obtained by the author in [13], and a different proof
is to be found in the dissertation of M. Lukkarinen [21].
A sketch of the proof now follows. Let
$$
{\bar L}_1(s) := \int_1^\infty|\z(\hf+iy)|^{2}{\roman e}^{-ys}
\d y\quad( \R s > 0). \leqno(5.2)
$$
Then we have by absolute convergence, taking $\s = \R s$ sufficiently large
and making the change of variable $xy = t$,
$$\eqalign{&
\int_0^\infty {\bar L}_1(x) x^{s-1}\d x = \int_0^\infty
\left(\int_1^\infty|\z(\hf + iy)|^{2}{\roman e}^{-yx}\d y\right)
x^{s-1}\d x\cr&
= \int_1^\infty|\z(\hf + iy)|^{2}\left(\int_0^\infty
x^{s-1}{\roman e}^{-xy}\d x\right)\d y\cr&
= \int_1^\infty|\z(\hf + iy)|^{2}y^{-s}\d y \int_0^\infty
{\roman e}^{-t}t^{s-1}\d t = {\Cal Z}_1(s)\G(s).\cr}\leqno(5.3)
$$
Further we have
$$\eqalign{&
\int_0^\infty {\bar L}_1(x) x^{s-1}\d x = \int_0^1{\bar L}_1(x) x^{s-1}\d x
+ \int_1^\infty{\bar L}_1(x) x^{s-1}\d x\cr&
= \int_1^\infty{\bar L}_1(1/x) x^{-1-s}\d x + A(s)\quad(\s>1),\cr}
$$
say, where $A(s)$ is an entire function. Since (see (1.1))
$$
{\bar L}_1(1/x) =  L_1(1/x) - \int_0^1|\z(\hf+iy)|^2{\roman e}^{-y/x}
\d y\qquad(x\ge1),
$$
it follows from (5.3) by analytic continuation that, for $\s>1$,
$$\eqalign{&
{\Cal Z}_1(s)\G(s) = \int_1^\infty L_1(1/x) x^{-1-s}\d x\cr&
- \int_1^\infty\bigl(\int_0^1|\z(\hf+iy)|^2{\roman e}^{-y/x}
\d y\bigr)x^{-1-s}\d x + A(s)\cr&
= I_1(s) - I_2(s) + A(s),\cr}\leqno(5.4)
$$
say. Clearly, for any integer $M\ge1$, we have
$$\eqalign{
I_2(s)&= \int_1^\infty\int_0^1|\z(\hf+iy)|^2
\left(\sum_{m=0}^M{(-1)^m\over m!}\left({y\over x}\right)^m +
O_M(x^{-M-1})\right)
\d y \,x^{-1-s}\d x\cr&
= \sum_{m=0}^M{(-1)^m\over m!}h_m\cdot{1\over m+s} + H_M(s),\cr}\leqno(5.5)
$$
say, where $H_M(s)$ is a regular function of $s$ for $\s > -M-1$, and
$h_m$ is a constant. Note that, for $\s = 1/T\,(T\to\infty)$
and any $N\ge0$, (1.5) gives
$$
L_1\left({1\over T}\right) =
\left(\log\left({T\over2\pi}\right) +
\gamma\right)\sum_{n=0}^N a_nT^{1-2n}
+ \sum_{n=0}^N b_nT^{-2n} + O_N(T^{-1-2N}\log T)
$$
with suitable $a_n, b_n\,(a_0 = 1)$. Inserting this formula in $I_1(s)$ in (5.4)
we have
$$\eqalign{
I_1(s) &= \int_1^\infty (\log {x\over2\pi}+\gamma)
\sum_{n=0}^Na_nx^{-2n-s}\d x + \int_1^\infty \sum_{n=0}^Nb_nx^{-1-2n-s}\d x
+ K_N(s)\cr&
= \sum_{n=0}^Na_n\left({1\over(2n+s-1)^2} +
{\gamma-\log2\pi\over 2n+s-1}\right)
+  K_N(s)\quad(\s > 1),\cr}\leqno(5.6)
$$
say, where $K_N(s)$ is regular for $\s > -2N$. Taking $M = 2N$
it follows from (5.4)--(5.6) that
$$\eqalign{
\Z_1(s)\G(s) &= \sum_{n=0}^Na_n\left({1\over(2n+s-1)^2} +
{\gamma-\log2\pi\over 2n+s-1}\right)\cr&
+ \sum_{m=0}^{2N}{(-1)^m\over m!}h_m\cdot{1\over m+s} + R_N(s),\cr}
\leqno(5.7)
$$
say, where $R_N(s)$ is a regular function of $s$ for $\s > -2N$. This holds
initially for $\s >1$, but by analytic continuation it holds for $\s > -2N$.
Since $N$ is arbitrary
and $\G(s)$ has no zeros, it follows that (5.7)
provides meromorphic continuation
of $\Z_1(s)$ to $\CC$. Taking into
account that $\G(s)$ has simple poles at $s = -m\;(m=0,1,2,\ldots\,)$
we obtain then the analytic
continuation of $\Z_1(s)$ to $\CC$, showing that besides $s=1$
the only poles of $\Z_1(s)$ can be
simple poles at $s = 1-2n$ for $n \in\NN$, as asserted by Theorem 4.
With more care the residues at these poles could be explicitly evaluated.
Finally using (5.7) and
$$
{1\over \G(s)} = 1 + \gamma(s-1) + \sum_{n=2}^\infty d_n(s-1)^n
$$
we obtain that the principal part of the Laurent expansion at $s=1$ is
given by (5.1).

\bigskip
Concerning the order of $\Z_1(s)$, we have (see M. Jutila [17])
$$
\Z_1(\s+it) \ll_\e t^{{5\over6}-\s+\e}\qquad(\hf \le \s \le 1,\;t\ge t_0).
\leqno(5.9)
$$
We also have the mean square bounds (see [16] for proof)
$$
\int_1^T|{\Cal Z}_1 (\s+it)|^2\d t \ll _\e T^{3-4\s+\e}\qquad(
0 \le \s \le{\txt{1\over2}}),\leqno(5.10)
$$
and
$$
\int_1^T|{\Cal Z}_1 (\s+it)|^2\d t \ll _\e
T^{2-2\s+\e}\qquad({\txt{1\over2}} \le \s \le 1).
\leqno(5.11)
$$
The bound in (5.11) is essentially best possible
since, for any given $\e > 0$,
$$
\int_1^T|\Z_k(\s + it)|^2\d t \;\gg_\e\; T^{2-2\s-\e}\qquad(k = 1,2;\;
\hf < \s < 1).\leqno(5.12)
$$
This assertion follows from
$$
\int_T^{2T}|\zt|^{2k}\d t \ll_\e
T^{2\s-1}\int_0^{T^{1+\e}}|\Z_k(\s+it)|^2\d t \quad(k = 1,2;\,\hf < \s < 1)
\leqno(5.13)
$$
and lower bounds for the integral on the left-hand side
(see [6, Chapter 9]).
The proof of (5.13) when $k = 2$  appeared in [11], and the
proof of the bound when $k=1$ is on similar lines.

\bigskip

\heading
{\bf 6. The  function $\Z_2  (s)$}
\endheading

The function ${\Cal Z}_2(s)$ has quite a different analytic behaviour
from the function $\Z_1(s)$. It
was introduced by Y. Motohashi [23]
(see also his monograph [24]). He  has shown   that ${\Cal Z}_2(s)$ has meromorphic
continuation over $\CC$. In the half-plane $\R s > 0$ it has
the following singularities: the pole $s = 1$ of order five,  simple
poles at $s = {1\over2} \pm i\k_j\,(\k_j =\sqrt{\lambda_j -
{1\over4}})$ and poles at $s =
\hf\rho$,  where $\rho$ denotes complex zeros of
$\zeta(s)$. The residue of ${\Cal Z}_2(s)$ at
$s = {1\over2} + i\k_h$ equals
$$ R(\k_h) := \sqrt{{\pi\over2}}{\Bigl(2^{-i\k_h}{\G({1\over4} -
\hf{i}\k_h)\over\G({1\over4} +
\hf{i}\k_h)}\Bigr)}^3\G(2i\k_h)\cosh(\pi\k_h)
\sum_{\k_j=\k_h}\alpha_j H_j^3(\txt{1\over2}),
$$
and the residue at $s = {1\over2} - i\k_h$ equals
$\overline{R(\k_h)}$.  Here as usual $\,\{\lambda_j = \k_j^2 +
{1\over4}\} \,\cup\, \{0\}\,$ is the discrete spectrum of the
non-Euclidean Laplacian acting on $SL(2,\ZZ)$-automorphic forms  and
$\a_j = |\rho_j(1)|^2(\cosh\pi\k_j)^{-1}$, where
$\rho_j(1)$ is the first Fourier coefficient of  the Maass wave form
corresponding to the eigenvalue $\lambda_j$ to which the Hecke
$L$-function
$H_j(s)$ is attached (see e.g., [24, Chapters 1--3]).
The principal part of ${\Cal Z}_2(s)$ has the form (this
may be compared with (5.1))
$$
\sum_{j=1}^5{A_j\over(s-1)^j},\leqno(6.1)
$$
where $A_5 = 12/\pi^2$, and the remaining $A_j$'s can be
evaluated explicitly by following the analysis in [23].
By Proposition 1 (see (4.4)) one obtains the leading term
in the asymptotic formula for the fourth moment of $|\zt|$,
although of course the coefficients in the full formula (i.e., (2.7)-(2.8)
for $k=2$) are well known explicitly (see the author's paper [8]).
The function ${\Cal Z}_2(s)$ was used to furnish several strong
results on $E_2(T)$ (see (2.7)), the error term in the asymptotic
formula for the fourth moment of $|\zt|$. Y. Motohashi [23]
used it to show that $E_2(T) = \Omega_\pm(T^{1/2})$,
which sharpens the earlier result of Ivi\'c-Motohashi (see [7])
that $E_2(T) = \Omega(T^{1/2})$.
The same authors (see e.g., [24]) have proved that
$$
E_2(T)  \ll T^{2/3}\log^{C_1}T\;(C_1>0),
\quad\int_0^TE_2(t)\d t \ll T^{3/2},\leqno(6.2)
$$
as well as the second bound in (2.12).  In [9] and [10] the
author has applied the theory of $\Z_2(s)$ to obtain the following
quantitative omega-results:
There exist constants $A, B > 0$ such that for $T \ge T_0 > 0$ every
interval $\,[T,\,AT]\,$ contains points $t_1, t_2, t_3, t_4$ such that
$$
E_2(t_1) > Bt_1^{1/2},\; E_2(t_2) < -Bt_2^{1/2},\;
\int_0^{t_3}E_2(t)\d t > Bt_3^{3/2},\;
\int_0^{t_4}E_2(t)\d t < -Bt_4^{3/2}.
$$
Moreover, we have (see [10])
$$
\int_0^TE_2^2(t)\d t \gg T^2,\leqno(6.3)
$$
which complements the upper bound in (6.2).

\medskip
As for the estimation of $\Z_2(s)$, we have (see [12])
$$
\Z_2(\s + it) \;\ll_\e\; t^{{4\over3}(1-\s)+\e}\qquad(\hf < \s \le 1;\,
t \ge t_0 > 0).\leqno(6.4)
$$
It was proved in [16] that
$$
\int_0^T|{\Cal Z}_2(\s + it)|^2\d t \ll_\e
T^\e\left(T + T^{2-2\s\over1-c}\right) \qquad(\hf < \s \le 1),\leqno(6.5)
$$
and we also have unconditionally
$$
\int_0^T|\Z_2(\s + it)|^2\d t \;\ll\;T^{10-8\s\over3}\log^CT
\qquad(\hf < \s \le 1,\,C > 0).\leqno(6.6)
$$
The constant $c$ appearing in (6.5) is defined by $
E_2(T) \ll_\e T^{c+\e},$ so that by (6.2) and (6.3) we have
$\hf \le c \le {2\over3}$. In (6.4)--(6.6) $\s$ is assumed
to be fixed, as $s = \s+it$ has to stay away from the $\hf$-line
where $\Z_2(s)$ has poles. Lastly, the author [14] proved that,
for ${5\over6} \le \s \le {5\over4}$ we have,
$$
\int_1^{T}|\Z_2(\s+it)|^2\d t \ll_\e T^{{15-12\s\over5}+\e}.\leqno(6.7)
$$
The lower limit of integration in (6.7) is unity, because of the pole
$s=1$ of $\Z_2(s)$. By taking $c = 2/3$ in (6.5) and using the convexity
of mean values (see [6, Lemma 8.3]) it follows that
$$
\int_1^T|\Z_2(\s+it)|^2\d t \;\ll_\e\; T^{{7-6\s\over2}+\e}\qquad(\hf < \s
\le \txt{5\over6}).\leqno(6.8)
$$
Note that (6.7) and (6.8) combined  provide the sharpest known bounds in the
whole range $\hf < \s \le {5\over6}$.

\medskip
Both pointwise and mean square estimates for $\Z_2(s)$ may be used
to estimate $E_2(T)$ and the eighth moment of $|\zt|$.
This connection is furnished by the following
result, proved by the author in [12].

\bigskip
THEOREM 5. {\it Suppose that, for some $\rho \ge 0$ and $r \ge 0$,}
$$
\Z_2(\s + it) \;\ll_\e\; |t|^{\rho+\e},\quad
\int_1^T|\Z_2(\s + it)|^2\d t \;\ll_\e\; T^{1+2r+\e}\quad(\hf < \s \le 1),
\leqno(6.9)
$$
{\it where $\s$ is fixed and $|t|\ge t_0 > 0$. Then we have}
$$
E_2(T) \;\ll_\e\; T^{{2\rho+1\over2\rho+2}+\e}, \quad
E_2(T) \;\ll_\e\; T^{{2r+1\over2r+2}+\e}\leqno(6.10)
$$
{\it and}
$$
\int_0^T|\zt|^8\d t \;\ll_\e\; T^{{4r+1\over2r+1}+\e}.\leqno(6.11)
$$
\medskip\no
Note that the conditions $\rho\ge0,\,r\ge 0$ must hold in view of (5.12).
Also note that from (6.6) with $\s = \hf + \e$ one can take in
(6.9) $r = \hf$, hence (6.10) gives
$E_2(T) \;\ll_\e\; T^{{2\over3}+\e},$
which is essentially the strongest known bound (see (6.2)). Thus
any improvement of the existing mean square bound for $\Z_2(s)$ at
$\s = \hf + \e$ would result in the bound for $E_2(T)$
 with the exponent strictly less than
2/3, which would be important.  Of course, if the first bound in (6.9)
holds with some $\rho$, then trivially the second bound will hold with
$r = \rho$, i.e. $r \le \rho$ has to hold. Observe that the
known value $r = \hf$ and (6.11) yield
$$
\int_0^T|\zt|^8\d t \;\ll_\e\; T^{\theta+\e}\leqno(6.12)
$$
with $\theta = 3/2$, which is, up to``$\e$", currently
the best known upper bound
(see [6, Chapter 8]) for the eighth moment, and any value $r < \hf$
in (6.9) would reduce the exponent $\theta = 3/2$ in (6.12). The connections
between upper bounds for the integral in (6.12) and mean
square estimates involving $\Z_2(s)$ and related functions are
also given in our last result. This is

\bigskip
THEOREM 6. {\it The eighth moment bound, namely} (6.12) {\it with
$\theta = 1$, is equivalent to the mean square bound
$$
\int_1^T|\Z_2(1+it)|^2\d t \;\ll_\e\; T^{\e},\leqno(6.13)
$$
and to
$$
\int_T^{2T}I^2(t,G)\d t \;\ll_\e\; T^{1+\e}\qquad(T^\e \le G \le T),\leqno(6.14)
$$
where}
$$
\quad I(T,G) := {1\over\sqrt{\pi}G}\int_{-\infty}^\infty
|\z(\hf + iT + iu)|^4{\roman e}^{-(u/G)^2}\d u.\leqno(6.15)
$$

\bigskip
{\bf Proof of Theorem 6}. We suppose first that (6.13) holds. Then
(6.12) with $\theta = 1$ follows from (5.13) with $k=2$.
Conversely, if (6.12) holds with $\theta = 1$,
note that we have, by [11, Lemma 4],
$$
\int_T^{2T}\left|\int_X^{2X}|\zx|^4x^{-s}\d x\right|^2\d t
\ll \int_X^{2X}|\zx|^8x^{1-2\s}\d x
\ll_\e X^{2-2\s+\e}\leqno(6.16)
$$
for $s = \s +it,\;\hf < \s \le 1$.
Similarly, using the Cauchy-Schwarz inequality for integrals
and the second bound in (2.12), it follows that
$$
\int_T^{2T}\left|\int_X^{2X}E_2(x)^2x^{-s}\d x\right|^2\d t
\ll_\e X^{1-2\s+\e}\qquad(s = \s + it,\;\s > \hf).\leqno(6.17)
$$
Combining (6.16) and (6.17) we obtain then,
similarly to the proof of (6.6) in [16],
$$
\int_1^T|\Z_2(\s+it)|^2\d t \;\ll_\e\; T^{4-4\s+\e}\qquad(\hf < \s\le 1),
$$
and we have (6.13). From (6.7) it follows that the integral in (6.13)
is unconditionally bounded by $T^{3/5+\e}$, and any improvement of the exponent
3/5 would also result in the improvement of the exponent
$\theta = 3/2$ in (6.12).

\medskip
Suppose again that (6.12) holds with $\theta = 1$. Then the left-hand
side of (6.14) is, for $T^\e \le G \le T$,
$$
\eqalign{&
\int_T^{2T}\left({1\over\sqrt{\pi}G}\int_{-\infty}^\infty
|\z(\hf + it + iu)|^4{\roman e}^{-(u/G)^2}\d u\right)^2\d t\cr&
\ll 1 + G^{-2}\int_T^{2T}\left(\int_{-G\log T}^{G\log T}
|\z(\hf + it + iu)|^4{\roman e}^{-(u/G)^2}\d u\right)^2\d t\cr&
\ll 1 + G^{-1}\int_{-G\log T}^{G\log T}
\left(\int_T^{2T}|\z(\hf + it + iu)|^8\d t\right)\d u\cr&
\ll_\e G^{-1}\int_{-G\log T}^{G\log T} T^{1+\e}\d u
\ll_\e  T^{1+\e},
\cr}
$$
as asserted. We remark that (6.14) is trivial when $T^{2/3} \le G \le T$
(see e.g., [7, Chapter 5]). Finally, if (6.14) holds, then we use
[15, Theorem 4], which in particular says that, for fixed $m\in\NN$,
$$
\int_T^{2T}I^m(t,G)\d t \ll_\e T^{1+\e}\quad(T^{\a_{m}+\e}\le G \le T,\,
 0 \le \a_{m} < 1)\leqno(6.18)
$$
implies that
$$
\int_0^T|\zt|^{4m}\d t \ll_\e T^{1+(m-1)\a_{m}+\e}.
$$
Using this result with $m = 2, \a_{2} = 0$, we obtain at once
(6.12) with $\theta = 1$. This completes the proof of Theorem 6.
So far (6.18) is known to hold unconditionally with $\a_2 = \hf$
(see [14]), which yields another proof of (6.12) with $\theta = 3/2$.

\bigskip
The significance of (6.14) is that for $I(T,G)$ in (6.15) an
explicit formula of Y. Motohashi (see [22] or [24]) exists. It
involves quantities from spectral theory, and thus the eighth
moment problem is directly connected to this theory via Theorem 6.
\Refs

\item{[1]} F.V. Atkinson, The mean value of the zeta-function on
the critical line, {\it Quart. J. Math. Oxford} {\bf 10}(1939), 122-128.

\item {[2]} F.V. Atkinson, The mean value of the zeta-function on
the critical line, {\it Proc. London Math. Soc.} {\bf 47}(1941), 174-200.

\item{[3]} J.B. Conrey, D.W, Farmer, J.P. Keating, M.O. Rubinstein and N.C.
Snaith, Integral moments of $L$-functions, {\it Proc. LMS} (to apear),
{\tt ArXiv:math.NT/0206018}.

\item{[4]} H. Delange, G\'en\'eralisation du th\'eor\`eme de Ikehara,
{\it Annales scien. E.N.S.} {\bf71}(1954), 213-242.

\item{[5]} A. Diaconu, D. Goldfeld and J. Hoffstein,
Multiple Dirichlet series and moments of zeta and $L$-functions,
{\it Compositio Math.} {\bf139}, No. 3. 297-360(2003).

\item {[6]} A. Ivi\'c,  The Riemann zeta-function, {\it John Wiley
and Sons}, New York, 1985.

\item {[7]} A. Ivi\'c,  Mean values of the Riemann zeta-function,
LN's {\bf 82}, {\it Tata Institute of Fundamental Research},
Bombay, 1991 (distr. by Springer Verlag, Berlin etc.).

\item{ [8]}  A. Ivi\'c,   On the fourth moment of the Riemann
zeta-function, {\it Publs. Inst. Math. (Belgrade)}
{ \bf57(71)}(1995), 101-110.

\item{[9]} A. Ivi\'c,  The Mellin transform and the Riemann
zeta-function,  {\it Proceedings of the Conference on Elementary and
Analytic Number Theory  (Vienna, July 18-20, 1996)},  Universit\"at
Wien \& Universit\"at f\"ur Bodenkultur, Eds. W.G. Nowak and J.
Schoi{\ss}engeier, Vienna 1996, 112-127.

\item{[10]} A. Ivi\'c, On the error term for the fourth moment of the
Riemann zeta-function, {\it J. London Math. Soc.}
{\bf60}(2)(1999), 21-32.

\item{ [11]}  A. Ivi\'c, On some conjectures and results
for the Riemann zeta-function
and Hecke series, {\it Acta Arith.}  {\bf109}(2001), 115-145.

\item{ [12]}  A. Ivi\'c, On the estimation of ${\Cal Z}_2(s)$,
in ``{\it Anal. Probab. Methods
Number Theory}" (eds. A. Dubickas et al.), TEV, Vilnius, 2002,  83-98.

\item{[13]} A. Ivi\'c, The Mellin transform of the square of Riemann's
zeta-function, International J. of Number Theory (to appear),
{\tt ArXiv:math.NT/0411404}.

\item{[14]} A. Ivi\'c, On the estimation of some Mellin transforms
connected with the fourth moment of $|\zt|$, subm. to the {\it Proc.
Conf. ``Elementare und Analytische Zahlentheorie"}, Mainz, May 2004
(ed. W. Schwarz), {\tt ArXiv:math.NT/0404524}.

\item{[15]} A. Ivi\'c, On moments of $|\zt|$ in short intervals,
subm. to the {\it Proc. Conf. ``Analytic Number Theory"} (Chennai,
December 2003), {\tt ArXiv:math.NT/0404289}.

\item{[16]} A. Ivi\'c, M. Jutila and Y. Motohashi, The Mellin transform of
power moments of the zeta-function, {\it Acta Arithmetica}
 {\bf95}(2000), 305-342.

\item {[17]}  M. Jutila, The Mellin transform of the square of Riemann's
zeta-function, {\it Periodica Math. Hung.} {\bf42}(2001), 179-190.

\item{[18]} J.P. Keating and  N.C. Snaith, Random Matrix Theory
and $L$-functions at $s=1/2$, {\it Comm. Math. Phys.} {\bf214}(2000), 57-89.

\item {[19]} H. Kober, Eine Mittelwertformel der Riemannschen Zetafunktion,
{\it Compositio Math.} {\bf3}(1936), 174-189.

\item {[20]} J. Korevaar, Tauberian Theory, {\it Grund. der math. Wissenschaften
Vol.} {\bf329}, Springer, Berlin etc., 2004.

\item{[21]} M. Lukkarinen, The Mellin transform of the square of Riemann's
zeta-function and Atkinson's formula, {\it Doctoral Dissertation}, University
of Turku, Turku, 2004, 69 pp.

\item{ [22]} Y. Motohashi,   An explicit formula for the fourth power
mean of the Riemann zeta-function, {\it Acta Math. }{\bf 170}(1993), 181-220.

\item {[23]} Y. Motohashi,  A relation between the Riemann zeta-function
and the hyperbolic Laplacian, {\it Annali Scuola Norm. Sup. Pisa, Cl. Sci. IV
ser.} {\bf 22}(1995), 299-313.

\item {[24]} Y. Motohashi,  Spectral theory of the Riemann
zeta-function, {\it Cambridge University Press}, Cambridge, 1997.

\item{[25]} M.A. Subhankulov, Tauberian theorems with remainder terms (Russian),
{\it Nauka}, Moscow, 1976.

\item{[26]} E.C. Titchmarsh, The theory of the Riemann zeta-function
(2nd ed.), Clarendon Press, Oxford, 1986.

\bigskip
\bigskip
\bigskip

Aleksandar Ivi\'c

Katedra Matematike RGF-a

Universitet u Beogradu, \DJ u\v sina 7

11000 Beograd, Serbia and Montenegro

\tt ivic\@rgf.bg.ac.yu, \enskip aivic\@matf.bg.ac.yu

\endRefs
\bye